\documentclass[12pt,a4paper]{article}
\usepackage{amsmath,amssymb}
\usepackage{graphicx}

\input{tcilatex}
\begin{document}

\title{Efficiency of the Multisection Method\thanks{%
Mathematics Subject Classifications: 65D99, 65Y20.}}
\date{{\small 20 February 2023}}
\author{J. S. C. Prentice\thanks{%
Senior Research Officer, Mathsophical Ltd., Johannesburg, South Africa}}
\maketitle

\begin{abstract}
We study the efficiency of the multisection method for univariate nonlinear
equations, relative to that for the well-known bisection \ method. We show
that there is a minimal effort algorithm that uses more sections than the
bisection method, although this optimal algorithm is problem dependent. The
number of sections required for optimality is determined by means of a
Lambert W function.
\end{abstract}

\section{Introduction}

Solving nonlinear equations computationally is an important part of applied
mathematics. The \textit{bisection method} \cite{B and F} is a simple yet
robust technique for doing so but, due to its linear convergence (in the
upper bound on its approximation error; see Appendix), it is usually
overlooked in comparison to faster methods, such as Newton's method. The 
\textit{mutlisection method} is an extension of the bisection method
(indeed, the bisection method is simply a special case of the multisection
method) designed to reduce the number of iterations, albeit at the expense
of increasing the number of function evaluations per iteration. It is
natural to wonder if, generally speaking, the multisection method is, in
some sense, more efficient than the bisection method. In this short paper we
investigate this notion by developing a suitable theoretical model, and
performing some appropriate numerical experiments.

\section{Relevant Concepts, Terminology and Notation}

The bisection method is an iterative method for solving the problem%
\begin{equation}
f\left( x\right) =0  \label{f(x)=0}
\end{equation}%
where $f\left( x\right) $ is real and continuous. If an interval $\left[ a,b%
\right] $ can be found on which $f\left( x\right) $ changes sign, then the
idea is to bisect $\left[ a,b\right] $ with a node $x_{1}$ (we will refer to
such an interval as the \textit{interval of relevance}). This node is taken
as an approximation to the root of \ref{f(x)=0} and, if $\left\vert f\left(
x_{1}\right) \right\vert $ is suitably small, $x_{1}$ is taken as the
solution to the problem. If $\left\vert f\left( x_{1}\right) \right\vert $
is not suitably small, then the signs of $f\left( x\right) $ at the three
nodes on the interval enable a new, smaller interval (either $\left[ a,x_{1}%
\right] $ or $\left[ x_{1},b\right] )$ to become the interval of relevance
in the next iteration. The process continues until $\left\vert f\left(
x_{i}\right) \right\vert $ is suitably small, or until the length of the
interval of relevance has been reduced to machine precision. This latter
condition is particularly useful, because it allows us to predict how many
iterations will be needed, at most. Note that each successive interval of
relevance is two times smaller than in the previous iteration.

The multisection method\ (also known as the $N$\textit{-section} method) is
similar to the bisection method, except that $N-1$ equispaced nodes are
introduced into the interval of relevance in each iteration $(N\in 
\mathbb{Z}
,$ $N\geqslant 2),$ thus creating $N$ equally sized subintervals (see, for
example, \cite{Mohd ali},where the cases $N=4$ and $N=6$ were considered).
If $\left\vert f\left( x\right) \right\vert $ is suitably small at any node,
the process is stopped and that node is taken as the solution. If not, a new
interval of relevance is identified and the process continues until the
interval of relevance has been reduced to machine precision. Each successive
interval of relevance is $N$ times smaller than in the previous iteration.
As in the bisection method, we can predict how many iterations will be
needed, at most. Note that the bisection method is simply the multisection
method with $N=2$.

In the next section, when determining the minimum number of iterations as a
function of $N,$ it will be necessary to solve an equation of the form%
\begin{equation}
ye^{y}=x,\text{ }x>0.  \label{y*exp(y)=x}
\end{equation}%
This equation has the solution%
\begin{equation*}
y=W\left( x\right)
\end{equation*}%
where $W\left( x\right) $ is the \textit{principal branch of the Lambert W
function} (sometimes denoted $W\left( 0,x\right) )$ \cite{Corless et al}\cite%
{Dence}. Additionally, $W\left( x\right) $ satisfies the identity%
\begin{equation}
e^{W\left( x\right) }=\frac{x}{W\left( x\right) }  \label{Lambert identity}
\end{equation}%
which will be useful later.

\section{Theory}

For the bisection method, the minimum number of iterations $M$ required to
reduce the interval of relevance to machine precision is given by%
\begin{equation*}
\frac{b-a}{2^{M}}\leqslant \mu \Rightarrow M\geqslant \frac{\ln \left( \frac{%
b-a}{\mu }\right) }{\ln 2}
\end{equation*}%
where we choose the smallest integer value of $M$ that satisfies this
inequality, and $\mu $ denotes machine precision (typically $2^{-52})$. The
analogous expression for the multisection method is%
\begin{equation*}
\frac{b-a}{N^{M}}\leqslant \mu \Rightarrow M\geqslant \frac{\ln \left( \frac{%
b-a}{\mu }\right) }{\ln N}.
\end{equation*}

Since the multisection method is iterative, a typical computational
implementation would involve a \texttt{for} or \texttt{while} loop.
Obviously, the number of such loops required is $M$. If we choose to measure
efficiency in terms of the total number of evaluations of $f\left( x\right)
, $ we would have $N-1$ evaluations of $f\left( x\right) $ per loop, giving
the total number of function evaluations $T_{f}$ as 
\begin{equation*}
T_{f}\left( N\right) =\frac{\left( N-1\right) \ln \left( \frac{b-a}{\mu }%
\right) }{\ln N}.
\end{equation*}%
However, this quantity does not capture the nuance of the algorithm. The
function $\left( N-1\right) /\ln N$ is monotonically increasing, so that $%
T_{f}\left( 2\right) $ is the least possible value of $T_{f}.$ This, in
turn, implies that the bisection method is more efficient than any other
possible multisection method. However, as we will see in the next section,
experimental results do not support this conclusion, so that $T_{f}$ is not
a good model for measuring the efficiency of the multisection method.

We prefer to measure efficiency in terms of total physical time $T_{t}.$ We
assume that the time taken to complete a single loop has the form $mN+c,$ so
that the total time taken to complete $M$ iterations is given by%
\begin{equation}
T_{t}\left( N\right) =\frac{\left( mN+c\right) \ln \left( \frac{b-a}{\mu }%
\right) }{\ln N}.  \label{Tt(N)}
\end{equation}%
Although this expression has the same form as that for $T_{f},$ we will see
that the constants $m$ and $c$ are vital in determining the presence of a
minimum (or not). We interpret $m$ as the cost (in time) of increasing $N$
by one, and $c$ as the cost (in time) of all parts of the loop that are not
dependent on $N$.

To find a minimum in $T_{t},$ we differentiate and equate to zero:%
\begin{align*}
\frac{dT_{t}}{dN}& =\ln \left( \frac{b-a}{\mu }\right) \left( \frac{m}{\ln N}%
-\frac{\left( mN+c\right) }{N\left( \ln N\right) ^{2}}\right) \\
\frac{dT_{t}}{dN}=0& \Rightarrow N\ln N=N+\frac{c}{m}.
\end{align*}

To solve this equation requires further manipulation. Using the notation $%
L\equiv \ln N,$ and the fact that $N=e^{\ln N},$ we find%
\begin{align*}
N\ln N& =N+\frac{c}{m} \\
\Rightarrow Le^{L}& =e^{L}+\frac{c}{m} \\
\Rightarrow e^{L}\left( L-1\right) & =\frac{c}{m} \\
\Rightarrow e^{L-1}\left( L-1\right) & =\frac{c}{m}e^{-1}.
\end{align*}%
The equation now has the form in (\ref{y*exp(y)=x}), with $y=L-1$ and $%
x=c/me.$ Hence, we have%
\begin{equation*}
L-1=W\left( \frac{c}{m}e^{-1}\right)
\end{equation*}%
so that%
\begin{align}
L& =\ln N=W\left( \frac{c}{m}e^{-1}\right) +1  \notag \\
\Rightarrow N& =e^{W\left( \frac{c}{m}e^{-1}\right) +1}=ee^{W\left( \frac{c}{%
m}e^{-1}\right) }  \notag \\
& =\frac{e\frac{c}{m}e^{-1}}{W\left( \frac{c}{m}e^{-1}\right) }  \notag \\
\Rightarrow N\left( R\right) & =\frac{R}{W\left( Re^{-1}\right) }\equiv N_{%
\text{\textit{min}}}  \label{expression for min N}
\end{align}%
where we made use of (\ref{Lambert identity}), and $R=c/m$. This is the
value of $N$ at which $T_{t}\left( N\right) $ has a minimum, and we have
emphasized the fact that it is dependent on the ratio $R$.

We now have%
\begin{align*}
T_{t}\left( N_{\text{\textit{min}}}\right) & =\frac{\left( mN_{\text{\textit{%
min}}}+c\right) \ln \left( \frac{b-a}{\mu }\right) }{\ln N_{\text{\textit{min%
}}}}=c\ln \left( \frac{b-a}{\mu }\right) \left( \frac{\frac{N_{\text{\textit{%
min}}}}{R}+1}{\ln N_{\text{\textit{min}}}}\right) \\
& =c\ln \left( \frac{b-a}{\mu }\right) \left( \frac{N_{\text{\textit{min}}}+R%
}{R\ln N_{\text{\textit{min}}}}\right)
\end{align*}%
which shows clearly that $T_{t}\left( N_{\text{\textit{min}}}\right) $ is
dependent not only on $R,$ but also on $a,b,c$ and $\mu .$ These latter four
parameters represent problem- and platform-specific aspects of the algorithm.

The ratio $T_{t}\left( N_{\text{\textit{min}}}\right) /T_{t}\left( 2\right) $%
, which measures efficiency at $N_{\text{\textit{min}}}$ relative to the
bisection method, and which we shall call $RelEff,$ is given by%
\begin{align}
RelEff\equiv \frac{T_{t}\left( N_{\text{\textit{min}}}\right) }{T_{t}\left(
2\right) }& =\frac{c\ln \left( \frac{b-a}{\mu }\right) \left( \frac{N_{\text{%
\textit{min}}}+R}{R\ln N_{\text{\textit{min}}}}\right) }{c\ln \left( \frac{%
b-a}{\mu }\right) \left( \frac{2+R}{R\ln 2}\right) }  \notag \\
& =\left( \frac{R\ln 2}{2+R}\right) \left( \frac{N_{\text{\textit{min}}}+R}{%
R\ln N_{\text{\textit{min}}}}\right) .  \label{RelEff}
\end{align}

We show $N_{\text{\textit{min}}}$ and $RelEff$ as functions of $R$ in
Figures 1 and 2. The minimum value of $N_{\text{\textit{min}}}$ is $2.73$.
The fact that it is not exactly $2$ (corresponding to the bisection method)
suggests that the bisection method is \textit{not} an optimal method in the
sense of our analysis. We see that $N_{\text{\textit{min}}}$ increases
monotonically, and we see that $RelEff$ achieves a minimum value of $\sim
0.1,$ for the values of $R$ considered here, suggesting the possibility of
significant improvements in efficiency relative to the bisection method.

\section{Numerical Calculations}

In this section, we report on numerical experiments pertaining to the theory
presented previously. We determine relevant quantities, such as $R,N_{\text{%
\textit{min}}}$ and $RelEff,$ for a selection of nonlinear equations, for $N$
in the range $N\in \left[ 2,250\right] .$ For each $N,$ we measure the time
per loop averaged over $1000$ loops. We summarize our results in Table 1.
The $r^{2}$ coefficient in Table 1 measures the goodness of fit for the
assumption that the time per loop can be assumed to be linear (i.e. $mN+c).$

\begin{center}
\begin{tabular}{|c|c|c|c|c|c|}
\hline
$f(x)$ & $[a,b]$ & $R$ $\left( =c/m\right) $ & $N_{\text{\textit{min}}}$ & $%
r^{2}$ & $RelEff$ \\ \hline
$\sin x-\cos x$ & $\left[ 0,\pi /2\right] $ & $2.73\times 10^{2}$ & $81$ & $%
0.963$ & $0.203$ \\ \hline
$e^{x}-2-x$ & $\left[ 1,4\right] $ & $2.17\times 10^{2}$ & $67$ & $0.980$ & $%
0.214$ \\ \hline
$x^{2}-8$ & $\left[ -5,-2\right] $ & $5.51\times 10^{2}$ & $140$ & $0.930$ & 
$0.175$ \\ \hline
$\frac{1}{10-x}-\frac{1}{4}$ & $\left[ -2,7\right] $ & $3.60\times 10^{2}$ & 
$100$ & $0.894$ & $0.191$ \\ \hline
$e^{-x}\sin x$ & $\left[ -8,-5\right] $ & $1.18\times 10^{2}$ & $43$ & $%
0.996 $ & $0.247$ \\ \hline
$\ln \left( x^{2}\right) -7$ & $\left[ 20,40\right] $ & $1.44\times 10^{2}$
& $49$ & $0.991$ & $0.236$ \\ \hline
\end{tabular}

Table 1: Results of numerical experiments.

\medskip
\end{center}

For each $f(x),$ we have determined $m$ and $c$ via a least-squares fit. The
coefficient $r^{2}$ indicates that these fits are very good. We have then
computed $N_{\text{\textit{min}}}$ and $RelEff$ using (\ref{expression for
min N}) and (\ref{RelEff}). It is clear that, for all cases, $N_{\text{%
\textit{min}}}$ is significantly greater than $2,$ and that $RelEff$ is $%
\sim 0.2$ $-$ some 80\% more efficient than the bisection method.

We provide some detail in Figures 3 and 4, where experimental data and
theoretical curves are shown for $f(x)=\sin x-\cos x.$ The other five cases
are similar. We note that $R$ is of the same order of magnitude for each
case and that, while $N_{\text{\textit{min}}}$ varies widely, $RelEff$ is
fairly consistent.

\section{Concluding Comments}

We have developed a model to describe the efficiency of the multisection
method. We measure efficiency in terms of the physical time required to
solve the given problem. Our model assumes that the time per iteration is
linear in $N,$ and numerical experiments seems to support this assumption.
The ratio of the parameters in this linear relationship determine the
optimal number of sections to be used, through the device of the Lambert W
function. It is clear that for all cases considered in our numerical study
that the optimal multisection method is significantly different to the
bisection method (i.e. $N_{\text{\textit{min}}}\gg 2$), and the resultant
improvement in efficiency is quite significant (about $80\%$). Of course,
our results are specific to our computational platform \cite{paltform}, but
we believe it is reasonable to assume that similar results would be obtained
on other systems.

\medskip

\section{Appendix}

Cheng and Lu \cite{Cheng} \ have shown, via a counter-example, that the
bisection method should, strictly speaking, \textit{not} be regarded as
being linearly convergent. While their argument is analytically correct, it
does seem to detract from the known property of the method that it will
yield a solution of acceptable accuracy within a finite number of
iterations. This property arises from the behaviour of the upper bound on
the approximation error of the method. Indeed, if $p$ is the root we seek,
and $p_{i}$ is the approximation after the $i$th iteration, we must have%
\begin{equation*}
0\leqslant \left\vert p_{i}-p\right\vert \leqslant \frac{b-a}{2^{i}}.
\end{equation*}%
The RHS of this expression is an upper bound on $\left\vert
p_{i}-p\right\vert ,$ which we will denote $B_{i},$ and we obviously have%
\begin{equation*}
B_{i+1}\leqslant \frac{1}{2}B_{i}
\end{equation*}%
for all $i\geqslant 0$ (with $B_{0}=b-a).$ \textit{The upper bound is
therefore linearly convergent}.

Moreover, there is always a value of $i$ such that 
\begin{equation*}
B_{i}\leqslant \varepsilon
\end{equation*}%
for any $0<\varepsilon \leqslant b-a$ so that, since $\left\vert
p_{i}-p\right\vert \leqslant B_{i},$ we have that $\left\vert
p_{i}-p\right\vert $ can be made as small as we choose. Furthermore, in a
finite precision computing environment, there exists a $z$ such that 
\begin{equation*}
B_{z}=\frac{b-a}{2^{z}}=0,
\end{equation*}%
so that for $z,$ $\left\vert p_{z}-p\right\vert =0$ and $p_{z}$ may be
considered exact in such an environment. On our platform, we find $z=1024.$

In \cite{Cheng}, an example is presented for which there is no integer $n$
such that the ratio%
\begin{equation*}
\frac{\left\vert p_{i+1}-p\right\vert }{\left\vert p_{i}-p\right\vert }
\end{equation*}%
decreases monotonically for all $i\geqslant n,$ and so the condition for
linear convergence is not satisfied. Nevertheless, the error bound \textit{%
does} converge (linearly) and the problem \textit{can} be solved in a finite
number of iterations (1024 iterations on our platform). In fact, we found $%
p=0.564468413605939$.

There is a tendency among practitioners to speak informally of the linear
convergence of the bisection method, simply because of this property of the
bound. More correctly, we should say that the bisection method is linearly
convergent in the upper bound of its approximation error. Nevertheless, the
bisection method is such that arbitrary accuracy \textit{can} be achieved in
a finite number of iterations and, in a finite precision environment, a 
\textit{de facto} exact result can be found.

It is natural to speculate that perhaps the notion of convergence in a
finite precision environment should be formally defined, but that is a topic
for another day.

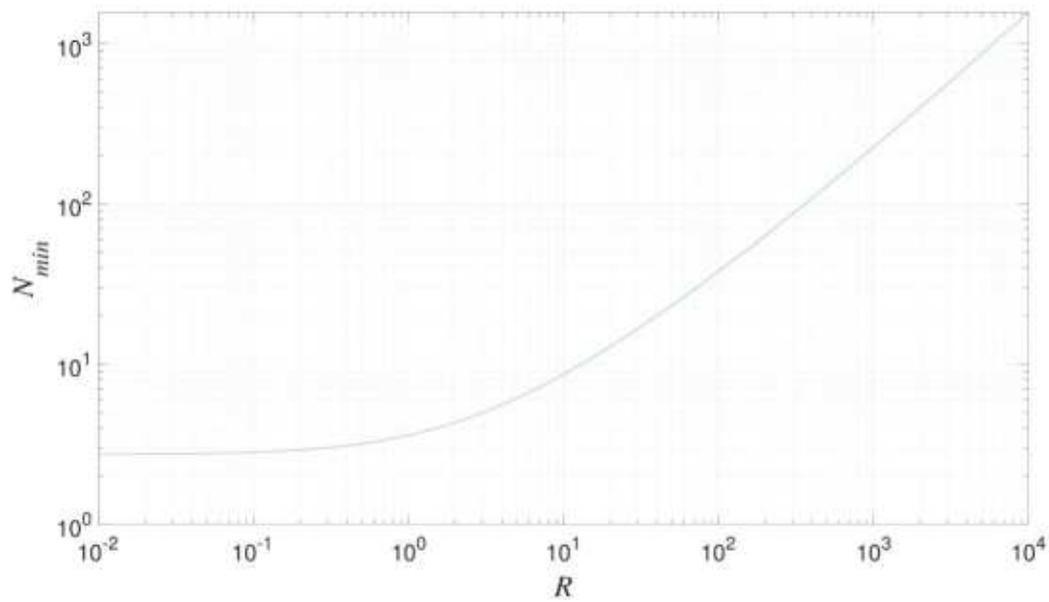
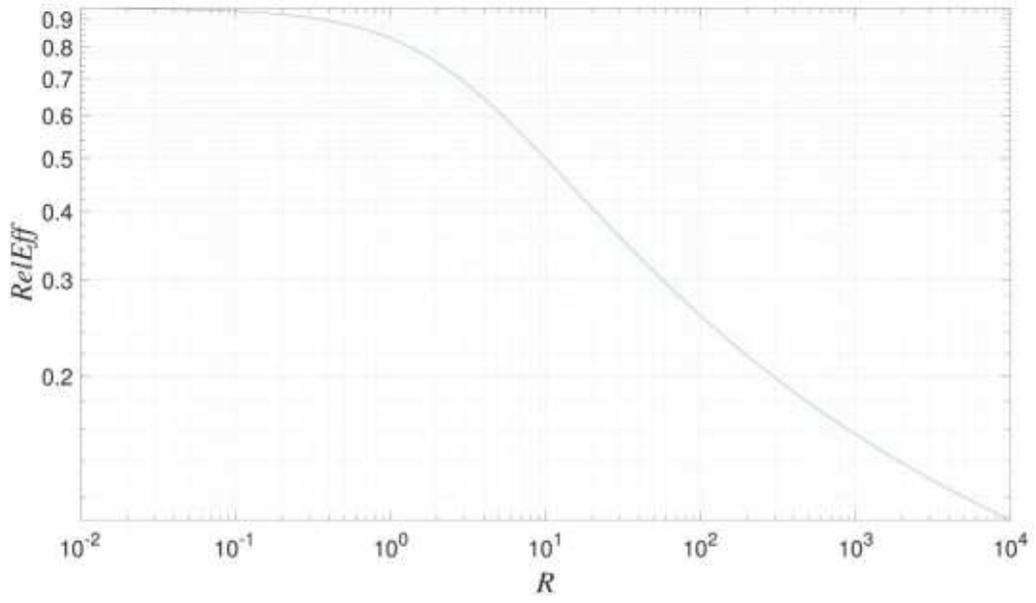
\end{document}